\documentclass{amsart}
\usepackage{amsfonts}
\usepackage{graphicx}

\newtheorem{theorem}{Theorem}[section]
\newtheorem*{theorem A}{Theorem A}
\newtheorem*{theorem B}{N\"olker's Theorem}

\theoremstyle{remark}

\theoremstyle{remark}

\theoremstyle{definition}

\newtheorem{example}{Example}[section]
\numberwithin{equation}{section}
\def\({\left ( }
\def\){\right )}
\def\<{\left < }
\def\>{\right >}

\input{tcilatex}

\begin{document}
\title{Constant curvature factorable surfaces in 3-dimensional isotropic space}
\author{Muhittin Evren Aydin}
\address{Department of Mathematics, Faculty of Science, Firat University, Elazig, 23200, Turkey}
\email{meaydin@firat.edu.tr}
\thanks{}

\begin{abstract}
In this paper, we study factorable surfaces in a
3-dimensional isotropic space. We classify such surfaces with constant isotropic Gaussian $(K)$ and mean curvature $(H)$. We provide a non-existence result related with the surfaces satisfying $\frac{H}{K}=const.$ Several examples are also illustrated.
\end{abstract}

\maketitle

\section{Introduction}

Let $\mathbb{E}^{3}$ be a 3-dimensional Euclidean space and $\left(
x,y,z\right) $ rectangular coordinates. A surface in $\mathbb{E}^{3}$ is
said to be \textit{factorable }(so-called \textit{homothetical}) if it is a
graph surface associated with $z=f\left( x\right) g\left( y\right) $ (see \cite{4,9}).

Constant Gaussian $(K)$ and mean curvature $(H)$ factorable surfaces in $%
\mathbb{E}^{3}$ were obtained in \cite{8,9,16}. As more general case, Zong
et al. \cite{17} defined that an \textit{affine factorable surface} in $%
\mathbb{E}^{3}$ is the graph of the function
\begin{equation*}
z=f\left( x\right) g\left(y+ax\right), \text{ } a\neq 0
\end{equation*} 
and classified ones with $K,H$ constants.

In a 3-dimensional Minkowski space $\mathbb{E}_{1}^{3},$ a surface is said to be 
\textit{factorable }if it can be expressed by one of the explicit forms:%
\begin{equation*}
\Phi _{1}:z=f\left( x\right) g\left( y\right), \text{ } \Phi _{2}:y=f\left( x\right) g\left( z\right), \text{ }\Phi _{3}:x=f\left( y\right) g\left( z\right) .
\end{equation*}

Six different classes of the factorable surfaces in $\mathbb{E}_{1}^{3}$ appear with respect to the causal characters of the directions (for
details, see \cite{10}). Such surfaces of $K,H$ constants in $\mathbb{E}%
_{1}^{3}$ were described in \cite{7,10,15}. 

Besides the Minkowskian space, a 3-dimensional isotropic space $\mathbb{I}%
^{3}$ provides two different types of the factorable surfaces. It is indeed
a product of the $xy-$plane and the isotropic $z-$direction with degenerate
metric (cf. \cite{5}). Due to the isotropic axes in $\mathbb{I}^{3}$ the factorable surface $\Phi _{1}$ distinctly behaves from others. We call it the \textit{factorable surface of type 1}. We refer to \cite{1}-\cite{3} for its details in $\mathbb{I}^{3}$.

The factorable surfaces $\Phi_{2}$, $\Phi _{3}$ in $\mathbb{I}^{3}$ are locally isometric
and, up to a sign, have same the second fundamental form. This means to have
same isotropic Gaussian $K$ and, up to a sign, mean curvature $H$. These surfaces are said to be of \textit{type 2}.

In this manner we are mainly interested with the factorable surfaces of type 2 in $\mathbb{I}^{3}$. We describe such surfaces in $%
\mathbb{I}^{3}$ with $K,H,H/K$ constants by the following results:

\begin{theorem}
A factorable surface of type 2 $\left( \Phi _{3}\right)$ in $\mathbb{I}^{3}$ has constant
isotropic mean curvature $H_{0}$ if and only if, up to suitable translations
and constants, one of the following occurs:

\begin{enumerate}
\item[(i)] If $\Phi _{3}$ is isotropic minimal, i.e. $H_{0}=0;$

\begin{enumerate}
\item[(i.1)] $\Phi _{3}$ is a non-isotropic plane,

\item[(i.2)] $x=y\tan \left( cz\right) $,

\item[(i.3)] $x=c\frac{z}{y}$.
\end{enumerate}

\item[(ii)] Otherwise $(H_{0}\neq 0),$ $x=\pm \sqrt{\frac{-z}{H_{0}}},$
\end{enumerate}
where $c$ is some nonzero constant.
\end{theorem}

\begin{theorem}
A factorable surface of type 2 $\left( \Phi _{3}\right)$ in $\mathbb{I}^{3}$ has constant
isotropic Gaussian curvature $K_{0}$ if and only if, up to suitable
translations and constants, one of the following holds:

\begin{enumerate}
\item[(i)] If $\Phi _{3}$ is isotropic flat, i.e. $K_{0}=0;$

\begin{enumerate}
\item[(i.1)] $x=c_{1}g\left( z\right), $  $g^{\prime }\neq 0$,

\item[(i.2)] $x=c_{1}e^{c_{2}y+c_{3}z}$,

\item[(i.3)] $x=c_{1}y^{c_{2}}z^{c_{3}},$ $c_{2}+c_{3}=1$.
\end{enumerate}

\item[(ii)] Otherwise $(K_{0}\neq 0);$ 
\begin{enumerate}
\item[(ii.1)] $K_{0}$ is negative and $x=\pm\frac{z}{ \sqrt{-K}y},$

\item[(ii.2)] 
\begin{equation*}
f\left( y\right) =\frac{c_{1}}{y} \text{ and } z=\pm \int \left( c_{2}g^{-1}-\frac{K_{0}}{c_{1}^{2}}\right)
^{-1/2}dg,
\end{equation*}

\end{enumerate}

\end{enumerate}
where $c_{1},c_{2},c_{3}$ are some nonzero constants.
\end{theorem}

\begin{theorem}
There does not exist a factorable surface of type 2 in $\mathbb{I}^{3}$ that satisfies $H/K=const. \neq 0$.
\end{theorem}

We remark that the results are also valid for the factorable surface
$\Phi _{2}$ in $\mathbb{I}^{3}$ by replacing $x$ with $y$ as well as taking $y=\pm\sqrt{z/H_{0}}$ in the last statement of Theorem 1.1.

\section{Preliminaries}

For detailed properties of isotropic spaces, see \cite{6}, \cite{11}-\cite%
{14}.

Let $P\left( \mathbb{R}^{3}\right) $ be a real 3-dimensional projective
space and $\omega $ a plane in $P\left( \mathbb{R}^{3}\right) .$ Then $%
P\left( \mathbb{R}^{3}\right) \backslash \omega $ becomes a real
3-dimensional affine space. Denote ($x_{0}:x_{1}:x_{2}:x_{3})\neq \left(
0:0:0:0\right) $ the projective coordinates in $P\left( \mathbb{R}%
^{3}\right) .$

A 3-dimensional \textit{isotropic space} $\mathbb{I}^{3}$ is an affine space
whose the absolute figure consists of a plane (absolute plane) $\omega $ and complex-conjugate straight lines (absolute lines) $l_{1},l_{2}$
in $\omega $. In coordinate form, $\omega $ is given by  $x_{0}=0$ and $l_{1},l_{2}$ by $x_{0}=x_{1}\pm ix_{2}=0.$ The absolute point, $\left(
0:0:0:1\right) ,$ is the intersection of the absolute lines.

For $x_{0}\neq 0,$ we have the affine coordinates by $x=\frac{x_{1}}{x_{0}},$
$y=\frac{x_{2}}{x_{0}},$ $z=\frac{x_{3}}{x_{0}}.$ The group of motions of $%
\mathbb{I}^{3}$ is given by
\begin{equation}
\left( x,y,z\right) \longmapsto \left( x^{\prime },y^{\prime },z^{\prime
}\right) :\left\{ 
\begin{array}{l}
x^{\prime }=a_{1}+x\cos \phi -y\sin \phi , \\ 
y^{\prime }=a_{2}+x\sin \phi +y\cos \phi , \\ 
z^{\prime }=a_{3}+a_{4}x+a_{5}y+z,%
\end{array}%
\right.  \tag{2.1}
\end{equation}%
where $a_{1},...,a_{5},\phi \in 
\mathbb{R}
.$

The \textit{isotropic metric} that is an invariant of $\left( 2.1\right)$ is induced by the absolute figure, namely $ds^{2}=dx^{2}+dy^{2}.$ One is degenerate along the
lines in $z-$direction and these lines are said to be \textit{isotropic.} A
plane is said to be \textit{isotropic} if it involves an isotropic line.
Otherwise it is called \textit{non-isotropic plane} or \textit{Euclidean
plane.}

We restrict our framework to regular surfaces whose the tangent plane at
each point is Euclidean, namely \textit{admissible} \textit{surfaces}.

Let $M$ be a regular admissible surface in $\mathbb{I}^{3}$ locally
parameterized by%
\begin{equation*}
r\left( u,v\right) =\left( x\left( u,v\right) ,y\left( u,v\right) ,z\left(
u,v\right) \right) 
\end{equation*}%
for a coordinate pair $(u,v)$. The components $E,F,G$ of
the first fundamental form of $M$ in $\mathbb{I}^{3}$ are computed by the
induced metric from $\mathbb{I}^{3}.$ The unit normal vector of $M$ is the
unit vector parallel to the $z-$direction.

The components of the second fundamental form $II$ of $M$ are given by%
\begin{equation*}
l=\frac{\det \left( r_{uu},r_{u},r_{v}\right) }{\sqrt{EG-F^{2}}},\text{ }m=%
\frac{\det \left( r_{uv},r_{u},r_{v}\right) }{\sqrt{EG-F^{2}}},\text{ }n=%
\frac{\det \left( r_{vv},r_{u},r_{v}\right) }{\sqrt{EG-F^{2}}}.
\end{equation*}
Accordingly, the \textit{isotropic
Gaussian} (or \textit{relative}) and \textit{mean curvature} of $M$ are respectively
defined by%
\begin{equation*}
K=\frac{ln-m^{2}}{EG-F^{2}},\text{ }H=\frac{En-2Fm+Gl}{2\left(
EG-F^{2}\right) }.
\end{equation*}%
A surface in $\mathbb{I}^{3}$ is said to be \textit{isotropic minimal} (\textit{%
isotropic flat}) if $H$ ($K$) vanishes.

\section{Proof of Theorem 1.1}
A factorable surface of type 2 in $\mathbb{I}^{3}$ can be locally
expressed by either%
\begin{equation*}
\Phi _{2}: r\left( x,z\right) =\left( x,f\left( x\right) g\left( z\right)
,z\right) 
\text{  or   }\Phi _{3}: r\left( y,z\right) =\left( f\left( y\right) g\left( z\right) ,y,z\right).
\end{equation*}%
All over this paper, all calculations shall be done for the surface $\Phi _{3}$. Its first fundamental form in $\mathbb{I}^{3}$ turns to 
\begin{equation*}
ds^{2}=\left( 1+\left( f^{\prime }g\right) ^{2}\right) dy^{2}+2\left(
fgf^{\prime }g^{\prime }\right) dydz+\left( fg^{\prime }\right) ^{2}dz^{2},%
\text{ }
\end{equation*}%
where $f^{\prime }=\frac{df}{dy}$, $g^{\prime }=\frac{dg}{dz}$.
Note that $g^{\prime }$ must be nonzero to obtain a regular admissible surface. By a calculation for the second fundamental form of $\Phi _{3}$ we have
\begin{equation*}
II=\left( \frac{f^{\prime \prime }g}{fg^{\prime }}\right) dy^{2}+2\left( 
\frac{f^{\prime }}{f}\right) dydz+\left( \frac{g^{\prime \prime }}{g^{\prime
}}\right) dz^{2}, \text{ } g^{\prime } \neq 0.
\end{equation*}%

Therefore the isotropic mean curvature $H$ of $\Phi _{3}$ becomes%
\begin{equation}
H=\frac{\left( \left( f^{\prime }g\right) ^{2}+1\right) g^{\prime \prime
}+g\left( g^{\prime }\right) ^{2}\left( ff^{\prime \prime }-2\left(
f^{\prime }\right) ^{2}\right) }{2f^{2}\left( g^{\prime }\right) ^{3}}. 
\tag{3.1}
\end{equation}%
Let us assume that $H=H_{0}=const.$ First we distinguish the case in which $\Phi _{3}$ is isotropic minimal:

\textbf{Case A. }$H_{0}=0.$ $\left( 3.1\right) $ reduces to%
\begin{equation}
\left( \left( f^{\prime }g\right) ^{2}+1\right) g^{\prime \prime }+g\left(
g^{\prime }\right) ^{2}\left( ff^{\prime \prime }-2\left( f^{\prime }\right)
^{2}\right) =0.  \tag{3.2}
\end{equation}

\textbf{Case A.1. }$f=f_{0}\neq 0\in \mathbb{R}$. $\left( 3.2\right) $
immediately implies $g=c_{1}z+c_{2},$ $c_{1},c_{2}\in \mathbb{R},$ and thus
we deduce that $\Phi _{3}$ is a non-isotropic plane. This gives the
statement (i.1) of Theorem 1.1.

\textbf{Case A.2. }$f=c_{1}y+c_{2},$ $c_{1},c_{2}\in \mathbb{R},$ $c_{1}\neq
0.$ $\left( 3.2\right) $ turns to%
\begin{equation*}
\frac{g^{\prime \prime }}{g^{\prime }}=\frac{2c_{1}^{2}gg^{\prime }}{%
1+\left( c_{1}g\right) ^{2}}.
\end{equation*}%
By solving this one, we obtain%
\begin{equation*}
g=\frac{1}{c_{1}}\tan \left( c_{2}z+c_{3} \right) ,\text{ 
}c_{2},c_{3}\in \mathbb{R},\text{ }c_{2}\neq 0,
\end{equation*}%
which proves the statement (i.2) of Theorem 1.1.

\textbf{Case A.3. }$f^{\prime \prime }\neq 0.$ By dividing $\left(
3.2\right) $ with $g\left( g^{\prime }\right) ^{2}$ one can be rewritten as%
\begin{equation}
\left( \left( f^{\prime }g\right) ^{2}+1\right) \frac{g^{\prime \prime }}{%
g\left( g^{\prime }\right) ^{2}}+ff^{\prime \prime }-2\left( f^{\prime
}\right) ^{2}=0.  \tag{3.3}
\end{equation}%
Taking partial derivative of $\left( 3.3\right) $ with respect to $z$ and
after dividing with $\left( f^{\prime }\right) ^{2},$ we get%
\begin{equation}
2\frac{g^{\prime \prime }}{g^{\prime }}+\left( \frac{1}{\left( f^{\prime
}\right) ^{2}}+g^{2}\right) \left( \frac{g^{\prime \prime }}{g\left(
g^{\prime }\right) ^{2}}\right) ^{\prime }=0.  \tag{3.4}
\end{equation}%
By taking partial derivative of $\left( 3.4\right) $ with respect to $y,$ we find $g^{\prime \prime }=c_{1}g\left(
g^{\prime }\right) ^{2},$ $c_{1}\in \mathbb{R}.$ We have two cases:

\textbf{Case A.3.1. }$c_{1}=0.$ $\left( 3.3\right) $ reduces to%
\begin{equation*}
ff^{\prime \prime }-2\left( f^{\prime }\right) ^{2}=0
\end{equation*}%
and by solving it we derive%
\begin{equation*}
f=-\frac{1}{c_{2}y+c_{3}},\text{ }c_{2},c_{3}\in \mathbb{R},\text{ }%
c_{2}\neq 0.
\end{equation*}%
This implies the statement (i.3) of Theorem 1.1.

\textbf{Case A.3.2. }$c_{1}\neq 0.$ It follows from $\left(
3.3\right) $ that%
\begin{equation}
c_{1}\left( \left( f^{\prime }g\right) ^{2}+1\right) +ff^{\prime \prime }-2\left(
f^{\prime }\right) ^{2}=0.  \tag{3.5}
\end{equation}%
Taking partial derivative of $\left( 3.5\right) $ with respect to $z$ yields 
$g^{\prime }=0$ which is not possible because of the regularity.

\textbf{Case B. }$H_{0}\neq 0.$ We have cases:

\textbf{Case B.1. }$f=f_{0}\neq 0\in \mathbb{R}$. Then $\left( 3.1\right) $
follows%
\begin{equation}
2H_{0}f_{0}^{2}=\frac{g^{\prime \prime }}{\left( g^{\prime }\right) ^{3}}, 
\tag{3.6.}
\end{equation}%
and solving it gives $g\left( z\right) =\pm \frac{1}{2H_{0}f_{0}^{2}}\sqrt{%
-4H_{0}f_{0}^{2}z+c_{1}}+c_{2},$ $c_{1},c_{2}\in \mathbb{R}$. This is the
proof of the statement (ii) of Theorem 1.1.

\textbf{Case B.2. }$f=c_{1}y+c_{2},$ $c_{1},c_{2}\in \mathbb{R},$ $c_{1}\neq
0.$ By considering this one in $\left( 3.1\right) $ we conclude%
\begin{equation}
2\left( c_{1}y+c_{2}\right) ^{2}H_{0}=\left( 1+c_{1}^{2}g^{2}\right) \frac{%
g^{\prime \prime }}{\left( g^{\prime }\right) ^{3}}-2c_{1}^{2}\frac{g}{%
g^{\prime }}.  \tag{3.7}
\end{equation}%
The left side in $\left( 3.7\right) $ is a function of $y$ while other side
is either a constant or a function $z.$ This is not possible.

\textbf{Case B.3. }$f^{\prime \prime }\neq 0.$ By multiplying both side of $\left(
3.1\right) $ with $2f^{2}\frac{g^{\prime }}{g}$ one can be rearranged as 
\begin{equation}
2H_{0}f^{2}\frac{g^{\prime }}{g}=\left( \left( f^{\prime }g\right)
^{2}+1\right) \frac{g^{\prime \prime }}{g\left( g^{\prime }\right) ^{2}}%
+ff^{\prime \prime }-2\left( f^{\prime }\right) ^{2}.  \tag{3.8}
\end{equation}%
Taking partial derivative of $\left( 3.8\right) $ with respect to $z$ and
after dividing with $\left( f^{\prime }\right) ^{2}$ yields%
\begin{equation}
2H_{0}\left( \frac{f}{f^{\prime }}\right) ^{2}\left( \frac{g^{\prime }}{g}%
\right) ^{\prime }=2\frac{g^{\prime \prime }}{g^{\prime }}+\left( g^{2}+%
\frac{1}{\left( f^{\prime }\right) ^{2}}\right) \left( \frac{g^{\prime
\prime }}{g\left( g^{\prime }\right) ^{2}}\right) ^{\prime }.  \tag{3.9}
\end{equation}%
After again taking partial derivative of $\left( 3.8\right) $ with respect
to $y$ we have%
\begin{equation}
2H_{0}\left( \left( \frac{f}{f^{\prime }}\right) ^{2}\right) ^{\prime
}\left( \frac{g^{\prime }}{g}\right) ^{\prime }=\left( \frac{1}{\left(
f^{\prime }\right) ^{2}}\right) ^{\prime }\left( \frac{g^{\prime \prime }}{%
g\left( g^{\prime }\right) ^{2}}\right) ^{\prime }.  \tag{3.10}
\end{equation}%
In order to solve (3.10) we have\ to consider several cases:

\textbf{Case B.3.1. }$f^{\prime }=c_{1}f,$ $c_{1}\in \mathbb{R},$ $c_{1}\neq
0.$ $\left( 3.10\right) $ leads to the following:

\textbf{Case B.3.1.1. }$g^{\prime \prime }=0,$ i.e, $g=c_{2}z+c_{3},$ $%
c_{2},c_{3}\in \mathbb{R},$ $c_{2}\neq 0.$ Then $\left( 3.8\right) $ reduces
to%
\begin{equation}
2H_{0}\frac{c_{2}}{c_{2}z+c_{3}}=-c_{1}^{2}, 
\tag{3.11}
\end{equation}%
which is a contradiction.

\textbf{Case B.3.1.2. }$g^{\prime \prime }=c_{2}g\left( g^{\prime }\right)
^{2},$ $c_{2}\in \mathbb{R},$ $c_{2}\neq 0.$ By dividing (3.8) with $f^{2}$ we get that%
\begin{equation*}
2H_{0}\frac{g^{\prime }}{g}=c_{1}^{2}c_{2}g^{2}+\frac{c_{2}}{f^{2}}%
-c_{1}^{2}
\end{equation*}%
and taking its partial derivative of $y$ gives the
contradiction $f^{\prime }=0.$

\textbf{Case B.3.2. }$f^{\prime }\neq c_{1}f,$ $c_{1}\in \mathbb{R}.$ If $%
g^{\prime }=c_{2}g,$ $c_{2}\in \mathbb{R},$ $c_{2}\neq 0$ in $\left(
3.10\right) $ then it follows%
\begin{equation*}
0=\left( \frac{1}{\left( f^{\prime }\right) ^{2}}\right) ^{\prime }\left( 
\frac{1}{g^{2}}\right) ^{\prime },
\end{equation*}%
which is not possible since $f^{\prime \prime }\neq 0$ and $g^{\prime }\neq
0.$ Hence $\left( 3.10\right) $ can be rewritten as%
\begin{equation}
2H_{0}\frac{\left( \left( \frac{f}{f^{\prime }}\right) ^{2}\right) ^{\prime }%
}{\left( \frac{1}{\left( f^{\prime }\right) ^{2}}\right) ^{\prime }}=\frac{%
\left( \frac{g^{\prime \prime }}{g\left( g^{\prime }\right) ^{2}}\right)
^{\prime }}{\left( \frac{g^{\prime }}{g}\right) ^{\prime }}.  \tag{3.12}
\end{equation}%
Both sides in $\left( 3.12\right) $ have to be a nonzero constant $c_{3}.$
Thereby $\left( 3.12\right) $ yields that%
\begin{equation}
\left( \frac{f}{f^{\prime }}\right) ^{2}=\frac{c_{3}}{\left( f^{\prime
}\right) ^{2}}+c_{4}  \tag{3.13}
\end{equation}%
and%
\begin{equation}
\frac{g^{\prime \prime }}{g\left( g^{\prime }\right) ^{2}}=2H_{0}c_{3}\frac{%
g^{\prime }}{g}+c_{5},  \tag{3.14}
\end{equation}%
where $c_{4},c_{5}\in \mathbb{R}.$ The fact that $f$  is a non-constant function leads to $c_{4} \neq 0$. $\left( 3.13\right) $ implies 
\begin{equation}
f^{\prime \prime }=\frac{1}{c_{4}}f.  \tag{3.15}
\end{equation}
Considering $\left( 3.13\right)-\left(3.15\right)$ in $\left( 3.8\right) $ gives%
\begin{equation}
2H_{0}f^{2}\frac{g^{\prime }}{g}=\left( \left( \frac{f^{2}-c_{3}}{c_{4}}%
\right) g^{2}+1\right) \left( 2H_{0}c_{3}\frac{g^{\prime }}{g}+c_{5}\right) -%
\frac{f^{2}}{c_{4}}+\frac{2c_{3}}{c_{4}}.  \tag{3.16}
\end{equation}%
By taking partial derivative of $\left( 3.16\right) $ with respect to $y,$
we find%
\begin{equation}
-2H_{0}\frac{g^{\prime }}{g}=\frac{c_{5}g^{2}-1}{c_{3}g^{2}-c_{4}}. 
\tag{3.17}
\end{equation}

\textbf{Case B.3.2.1. }$c_{5}=0.$ Then $\left( 3.14\right) $ follows%
\begin{equation}
g^{\prime }=\frac{-1}{2H_{0}c_{3}g+c_{6}},\text{ }c_{6}\in \mathbb{R}. 
\tag{3.18}
\end{equation}%
Substituting $\left( 3.18\right) $ in $\left( 3.17\right) $ leads to%
\begin{equation*}
\frac{2H_{0}}{2H_{0}c_{3}g^{2}+c_{6}g}=\frac{-1}{c_{3}g^{2}-c_{4}}
\end{equation*}%
or the following polynomial equation on $g$:%
\begin{equation}
4H_{0}c_{3}g^{2}+c_{6}g-2H_{0}=0.  \tag{3.19}
\end{equation}%
Since the coefficients $H_{0}$ and $c_{3}$ are nonzero, we obtain a
contradiction.

\textbf{Case B.3.2.2. }$c_{5}\neq 0.$ Since $\frac{g^{\prime }}{g}$ is not
constant, $\left( 3.17\right) $ immediately implies $c_{4}\neq \frac{c_{3}}{%
c_{5}}$. Substituting $\left( 3.17\right) $ into $\left( 3.14\right) $ gives%
\begin{equation*}
\frac{g^{\prime \prime }}{g^{\prime }}=\left( c_{3}-c_{4}c_{5}\right) \frac{%
gg^{\prime }}{c_{3}g^{2}-c_{4}}. 
\end{equation*}%
or
\begin{equation}
g^{\prime }=c_{7}\left( c_{3}g^{2}-c_{4}\right) ^{\frac{c_{3}-c_{4}c_{5}}{%
2c_{3}}}, \text{ } c_{7} \neq 0.  \tag{3.20}
\end{equation}%
By considering $\left( 3.20\right) $ in $\left( 3.17\right) $ we deduce%
\begin{equation}
-2H_{0}c_{6}\left( c_{3}g^{2}-c_{4}\right) ^{\frac{3c_{3}-c_{4}c_{5}}{2c_{3}}%
}=c_{5}g^{3}-g.  \tag{3.21}
\end{equation}%
This leads to a contradiction since the terms $g$ of different
degrees appears in $\left( 3.21\right) .$

\section{Proof of Theorem 1.2}

By a calculation for a factorable graph of type 2 in 
$\mathbb{I}^{3}$, the isotropic Gaussian curvature turns to%
\begin{equation}
K=\frac{fgf^{\prime \prime }g^{\prime \prime }-\left( f^{\prime }g^{\prime
}\right) ^{2}}{\left( fg^{\prime }\right) ^{4}}.  \tag{4.1}
\end{equation}%
Let us assume that $K=K_{0}=const.$ We have cases:

\textbf{Case A.} $K_{0}=0.$ $\left( 4.1\right) $ reduces to%
\begin{equation}
fgf^{\prime \prime }g^{\prime \prime }-\left( f^{\prime }g^{\prime }\right)
^{2}=0.  \tag{4.2}
\end{equation}
$f$ or $g$ constants are solutions for $\left( 4.2\right) $ and by regularity we have
the statement (i.1) of Theorem 1.2. Suppose that $f,g$ are non-constants. Then $%
\left( 4.2\right) $ yields $f^{\prime \prime }g^{\prime \prime }\neq 0.$
Thereby $\left( 4.2\right) $ can be arranged as%
\begin{equation}
\frac{ff^{\prime \prime }}{\left( f^{\prime }\right) ^{2}}=\frac{\left(
g^{\prime }\right) ^{2}}{gg^{\prime \prime }}.  \tag{4.3}
\end{equation}%
Both sides of $\left( 4.3\right) $ are equal to same nonzero constant, namely%
\begin{equation}
ff^{\prime \prime }-c_{1}\left( f^{\prime }\right) ^{2}=0 \text{ and }%
gg^{\prime \prime }-\frac{1}{c_{1}}\left( g^{\prime }\right) ^{2}=0. 
\tag{4.4}
\end{equation}%
If $c_{1}=1$ in $\left( 4.4\right) ,$ then by solving it
we obtain%
\begin{equation*}
f\left( y\right) =c_{2}e^{c_{3}y}\text{ and }g\left( z\right)
=c_{4}e^{c_{5}z},\text{ }c_{2},...,c_{5}\in \mathbb{R}.
\end{equation*}%
This gives the statement (i.2) of Theorem 1.2. Otherwise, i.e. $c_{1}\neq 1$
in $\left( 4.4\right) ,$ we derive%
\begin{equation*}
f\left( y\right) =\left( \left( 1-c_{1}\right) \left( c_{6}y+c_{7}\right)
\right) ^{\frac{1}{1-c_{1}}}\text{ and }g\left( z\right) =\left( \left( 
\frac{c_{1}-1}{c_{1}}\right) \left( c_{8}z+c_{9}\right) \right) ^{\frac{c_{1}%
}{c_{1}-1}},
\end{equation*}%
where $c_{6},...,c_{9}\in \mathbb{R}.$ This completes the proof of the
statement (i) of Theorem 1.2.

\textbf{Case B.} $K_{0}\neq 0.$ $\left( 4.1\right) $ can be rewritten as%
\begin{equation}
K_{0}\left( g^{\prime }\right) ^{2}=\frac{f^{\prime \prime }}{f^{3}}\left( 
\frac{gg^{\prime \prime }}{\left( g^{\prime }\right) ^{2}}\right) -\left( 
\frac{f^{\prime }}{f^{2}}\right) ^{2}.
  \tag{4.5}
\end{equation}%
Taking parital derivative of $\left( 4.5\right) $ with respect to $z$ leads
to%
\begin{equation}
2K_{0}g^{\prime }g^{\prime \prime }=\frac{f^{\prime \prime }}{f^{3}}\left( 
\frac{gg^{\prime \prime }}{\left( g^{\prime }\right) ^{2}}\right) ^{\prime }.
\tag{4.6}
\end{equation}%
We have several cases for $\left( 4.6\right) $:

\textbf{Case B.1.} $g^{\prime \prime }=0,$ $g\left( z\right) =c_{1}z+c_{2},$ 
$c_{1},c_{2}\in \mathbb{R}.$ Hence from $%
\left( 4.5\right) $ we deduce%
\begin{equation*}
K_{0}\left( c_{1}\right) ^{2}=-\left( \frac{f^{\prime }}{f^{2}}\right) ^{2},
\end{equation*}%
which implies that $K_{0}$ is negative and 
\begin{equation*}
f\left( y\right) =\frac{1}{\pm c_{1}\sqrt{-K_{0}}y+c_{3}}.
\end{equation*}%
This proves the statement (ii.1) of Theorem 1.2.

\textbf{Case B.2.} $g^{\prime \prime }\neq 0$. $\left( 4.6\right) $
immediately implies $f^{\prime \prime }\neq 0$. Then taking parital
derivative of $\left( 4.6\right) $ with respect to $y$ gives%
\begin{equation}
0=\left( \frac{f^{\prime \prime }}{f^{3}}\right) ^{\prime }\left( \frac{%
gg^{\prime \prime }}{\left( g^{\prime }\right) ^{2}}\right) ^{\prime }, 
\tag{4.7}
\end{equation}%
or
\begin{equation}
f^{\prime \prime }=c_{1}f^{3}, \text{ } c_{1} \in \mathbb{R}. \tag{4.8}
\end{equation}
By considering (4.8) in $\left( 4.5\right) $ we get%
\begin{equation}
K_{0}\left( g^{\prime }\right) ^{2}=c_{1} \frac{%
gg^{\prime \prime }}{\left( g^{\prime }\right) ^{2}}-\left( \frac{f^{\prime }}{%
f^{2}}\right) ^{2}.  \tag{4.9}
\end{equation}%
Taking partial derivative of $\left( 4.9\right) $ with respect to $y$ leads to 
\begin{equation}
f^{\prime }=c_{2}f^{2}, \text{ } c_{2} \in \mathbb{R}. \tag{4.10}
\end{equation}
It follows from (4.8) and (4.10) that $c_{1}=2c_{2}^{2}$ and 
\begin{equation*}
f\left( y\right) =-\frac{1}{c_{2}y+c_{3}}
\end{equation*}%
for some constant $c_{3}.$ By substituting (4.8) and (4.10) into (4.5), we conclude
\begin{equation}
\frac{K_{0}}{c_{2}^{2}}r^{3}+r=2gr^{\prime }, \tag{4.11}
\end{equation}%
where $r=g^{\prime }$ and $r^{\prime }=\frac{dr}{dg}=\frac{g^{\prime \prime }%
}{g^{\prime }}.$ After solving $\left( 4.11\right) ,$ we obtain%
\begin{equation*}
r=\pm \left( c_{4}^{2}g^{-1}-\frac{K_{0}}{c_{2}^{2}}\right) ^{-1/2}, \text{ } c_{4} \in \mathbb{R}, \text{ } c_{4} \neq 0,
\end{equation*}
or%
\begin{equation*}
z=\pm \int \left( c_{4}^2g^{-1}-\frac{K_{0}}{c_{2}^{2}}\right)
^{1/2}dg,
\end{equation*}%
which proves the statement (ii.2) of Theorem 1.2.

\section{Proof of Theorem 1.3}

Assume that a factorable surface of type 2 in $\mathbb{I}^{3}$
fulfills the condition $H+ \lambda K=0,$ $\lambda HK\neq 0.$ Then $\left( 3.1\right) $ and $%
\left( 4.1\right) $ give rise to%
\begin{equation}
f^{2}\left( \left( f^{\prime }g\right) ^{2}+1\right) g^{\prime }g^{\prime
\prime }+f^{2}g\left( g^{\prime }\right) ^{3}\left( ff^{\prime \prime
}-2\left( f^{\prime }\right) ^{2}\right) +2\lambda \left( ff^{\prime \prime
}gg^{\prime \prime }-\left( f^{\prime }g^{\prime }\right) ^{2}\right) =0. 
\tag{5.1}
\end{equation}%
Since $K\neq 0,$ $\left( 5.1\right) $ can be divided by $\left( ff^{\prime
}\right) ^{2}$ as follows:%
\begin{equation}
\left( g^{2}+\frac{1}{\left( f^{\prime }\right) ^{2}}\right) g^{\prime
}g^{\prime \prime }+g\left( g^{\prime }\right) ^{3}\left( \frac{ff^{\prime
\prime }}{\left( f^{\prime }\right) ^{2}}-2\right) +2\lambda \frac{f^{\prime \prime
}}{f\left( f^{\prime }\right) ^{2}}gg^{\prime \prime }-\frac{2\lambda }{f^{2}}%
\left( g^{\prime }\right) ^{2}=0.  \tag{5.2}
\end{equation}%
In order to solve $\left( 5.2\right) $ we have to distinguish several cases:

\textbf{Case A.} $g=c_{1}z+c_{2}$, $c_{1},c_{2} \in \mathbb{R}$, $c_{1}\neq 0.$ $\left(5.2\right) $
reduces to%
\begin{equation}
c_{1}\left( c_{1}z+c_{2}\right) \left( \frac{ff^{\prime \prime }}{\left(
f^{\prime }\right) ^{2}}-2\right) -\frac{2\lambda }{f^{2}}=0.  \tag{5.3}
\end{equation}%
Taking partial derivative of $\left( 5.3\right) $ with respect to $z$ gives $%
ff^{\prime \prime }=2\left( f^{\prime }\right) ^{2}.$ Considering it into $%
\left( 5.3\right) $ yields the contradiction $\lambda = 0$.

\textbf{Case B.} $g^{\prime \prime }\neq 0.$ By dividing $\left(5.2\right) $
with the product $g^{\prime }g^{\prime \prime },$ we get%
\begin{equation}
g^{2}+\frac{1}{\left( f^{\prime }\right) ^{2}}+\frac{g\left( g^{\prime
}\right) ^{2}}{g^{\prime \prime }}\left( \frac{ff^{\prime \prime }}{\left(
f^{\prime }\right) ^{2}}-2\right) +2\lambda \frac{f^{\prime \prime }g}{f\left(
f^{\prime }\right) ^{2}g^{\prime }}-\frac{2\lambda g^{\prime }}{f^{2}g^{\prime
\prime }}=0.  \tag{5.4}
\end{equation}%
Put $p=f^{\prime },$ $p^{\prime }=\frac{dp}{df}=\frac{f^{\prime \prime }}{%
f^{\prime }},$ $r=g^{\prime }$ and $r^{\prime }=\frac{dr}{dg}=\frac{%
g^{\prime \prime }}{g^{\prime }}$ in (5.4). Thus taking partial derivatives of $\left(
5.4\right) $ with respect to $f$ and $g$ implies%
\begin{equation}
\left( \frac{gr}{r^{\prime }}\right) ^{\prime }\left( \frac{fp^{\prime }}{p}%
\right) ^{\prime }+2\lambda \left( \frac{p^{\prime }}{fp}\right) ^{\prime }\left( 
\frac{g}{r}\right) ^{\prime }-2\lambda \left( \frac{1}{f^{2}}\right) ^{\prime
}\left( \frac{1}{r^{\prime }}\right) ^{\prime }=0.  \tag{5.5}
\end{equation}%
Since $\left( 1/f^{2}\right) ^{\prime }\neq 0,$ we can rewrite $\left(
5.5\right) $ as%
\begin{equation}
\left( \frac{gr}{r^{\prime }}\right) ^{\prime }\frac{\left( fp^{\prime
}/p\right) ^{\prime }}{\left( 1/f^{2}\right) ^{\prime }}+2\lambda\left( \frac{g}{r}%
\right) ^{\prime }\frac{\left( p^{\prime }/fp\right) ^{\prime }}{\left(
1/f^{2}\right) ^{\prime }}-2\lambda \left( \frac{1}{r^{\prime }}\right) ^{\prime
}=0.  \tag{5.6}
\end{equation}%
Taking derivative of $\left(
5.6\right) $ with respect to $f$ leads to%
\begin{equation}
\left( \frac{gr}{r^{\prime }}\right) ^{\prime }\left( \frac{\left(
fp^{\prime }/p\right) ^{\prime }}{\left( 1/f^{2}\right) ^{\prime }}\right)
^{\prime }+2\lambda \left( \frac{g}{r}\right) ^{\prime }\left( \frac{\left(
p^{\prime }/fp\right) ^{\prime }}{\left( 1/f^{2}\right) ^{\prime }}\right)
^{\prime }=0.  \tag{5.7}
\end{equation}%
We have some cases to solve (5.7):

\textbf{Case B.1.} $p^{\prime }=0$, i.e. $f\left( y\right) =c_{1}y+c_{2,}$ $c_{1},c_{2} \in \mathbb{R}$, $c_{1}\neq 0.$ Considering it into $\left( 5.4\right) $ gives%
\begin{equation}
g^{2}+\frac{1}{c_{1}^{2}}-2\frac{g\left( g^{\prime }\right) ^{2}}{g^{\prime
\prime }}-\frac{2\lambda g^{\prime }}{f^{2}g^{\prime \prime }}=0  \tag{5.8}
\end{equation}%
and taking partial derivative of $\left( 5.8\right) $ with respect to $y$
implies that $f^{\prime }$ or $g^{\prime }$ vanish however both situations
are not possible.

\textbf{Case B.2} In $\left( 5.7\right) $ assume that $p^{\prime } \neq 0$ and%
\begin{equation}
\left( \frac{\left( fp^{\prime }/p\right) ^{\prime }}{\left( 1/f^{2}\right)
^{\prime }}\right) ^{\prime }=\left( \frac{\left( p^{\prime }/fp\right)
^{\prime }}{\left( 1/f^{2}\right) ^{\prime }}\right) ^{\prime }=0.  \tag{5.9}
\end{equation}%
This one follows%
\begin{equation}
\frac{fp^{\prime }}{p}=\frac{c_{1}}{f^{2}}+c_{2}\text{ and }\frac{p^{\prime }%
}{fp}=\frac{c_{3}}{f^{2}}+c_{4},\text{ }c_{1},...,c_{4}\in \mathbb{R}. 
\tag{5.10}
\end{equation}%
Both equalities in $\left( 5.10\right) $ imply%
\begin{equation}
\frac{p^{\prime }}{p}=\frac{c_{2}}{f^{2}},\text{ }c_{2}\neq 0,\text{ }%
c_{2}=c_{3}.  \tag{5.11}
\end{equation}%
Considering $\left( 5.11\right) $ in the first or second equality of $\left(
5.9 \right) $ leads to a contradiction.

\textbf{Case B.3} $\left( g/r\right) ^{\prime }=0.$ This implies $%
g^{\prime }=c_{1}g,$ namely $g=c_{2}e^{c_{1}z},$ $c_{1},c_{2} \in \mathbb{R}$, $c_{1}c_{2}\neq 0.$
Substituting it into $\left( 5.4\right) $ yields%
\begin{equation*}
f\left( y\right) =c_{3}e^{c_{4}y},  
\end{equation*}%
which is a contradiction since $K\neq 0.$

\textbf{Case B.4.} $\left( gr/r^{\prime }\right) ^{\prime }=0$
in $\left( 5.7\right) .$ Then $r^{\prime }=c_{1}gr,$ $c_{1} \in \mathbb{R}$, $c_{1}\neq 0$, and $%
\left( 5.6\right) $ reduces to%
\begin{equation}
\left( \frac{g}{r}\right) ^{\prime }\frac{\left( p^{\prime }/fp\right)
^{\prime }}{\left( 1/f^{2}\right) ^{\prime }}-\left( \frac{1}{r^{\prime }}%
\right) ^{\prime }=0,  \tag{5.12}
\end{equation}%
and taking partial derivative of $\left( 5.12\right) $ with respect to $g$
yields%
\begin{equation}
\frac{\left( p^{\prime }/fp\right) ^{\prime }}{\left( 1/f^{2}\right)
^{\prime }}=c_{2}\neq 0.  \tag{5.13}
\end{equation}%
Substituting $\left(5.13\right) $ into $\left( 5.12\right) $ gives%
\begin{equation}
r=c_{3}g-\frac{1}{c_{4}g},\text{ }c_{3},c_{4} \in \mathbb{R}, \text{ }c_{3}c_{4}\neq 0.  \tag{5.14}
\end{equation}%
By taking derivative of (5.14) with respect to $g$ and comparing
with $r^{\prime }=c_{1}gr,$ we deduce the following polynomial equation on $g$:
\begin{equation*}
c_{1}c_{3}g^{4}-\left( \frac{c_{1}}{c_{4}}+c_{3}\right) g^{2}-\frac{1}{c_{2}}%
=0.
\end{equation*}
This gives a contradiction.

\textbf{Case B.5.} $\left( gr/r^{\prime }\right) ^{\prime } \neq 0$
in $\left( 5.8\right) .$ Then (5.8) can be rewritten as 
\begin{equation*}
\underset{G\left( g\right) }{\underbrace{\left( \left( \frac{g}{r}\right)
^{\prime }\right) ^{-1}\left( \frac{gr}{r^{\prime }}\right) ^{\prime }}}%
+2\lambda \underset{F\left( f\right) }{\underbrace{\left( \frac{\left(
p^{\prime }/fp\right) ^{\prime }}{\left( 1/f^{2}\right) ^{\prime }}\right)
^{\prime }\left( \left( \frac{\left( fp^{\prime }/p\right) ^{\prime }}{%
\left( 1/f^{2}\right) ^{\prime }}\right) ^{\prime }\right) ^{-1}}}=0,
\end{equation*}
which implies that $G(g)=c_{1},$ $F(f)=-c_{1}/2\lambda$, $c_{1}\in \mathbb{%
R}$, $c_{1}\neq 0$. Thus we have
\begin{equation}
\frac{gr}{r^{\prime }}=c_{1}\frac{g}{r}+c_{2},\text{ }c_{2}\in \mathbb{%
R}.  \tag{5.15}
\end{equation}%
Substituting (5.15) in (5.5) follows%
\begin{equation}
c_{1}\left( \frac{fp^{\prime }}{p}\right) ^{\prime }+2\lambda \left( \frac{%
p^{\prime }}{fp}\right) ^{\prime }-2\lambda \left( \frac{1}{f^{2}}\right) 
\frac{\left( 1/r^{\prime }\right) ^{\prime }}{\left( g/r\right) ^{\prime }}%
=0.  \tag{5.16}
\end{equation}%
By taking partial derivative of (5.16) with respect to $g$, we find%
\begin{equation}
\frac{g}{r}=\frac{c_{3}}{r^{\prime }}+c_{4},\text{ }c_{3},c_{4}\in \mathbb{R}%
,\text{ }c_{4}\neq 0.  \tag{5.17}
\end{equation}%
Substituting (5.15) and (5.17) into (5.4) gives 
\begin{equation}
\frac{g}{r}=c_{5}g^{2}+c_{6},\text{ }c_{5},c_{6}\in \mathbb{R},\text{ }%
c_{5}\neq 0.  \tag{5.18}
\end{equation}%
$\left( 5.15\right) ,$ $\left( 5.17\right) $ and $\left( 5.18\right) $ imply
the following polynomial equation:%
\begin{equation*}
\left( c_{1}c_{5}-c_{3}c_{5}^{2}\right) g^{4}+\left(
c_{1}c_{6}-2c_{3}c_{5}c_{6}-\left( c_{2}-c_{4}\right) c_{5}\right)
g^{2}+\left( c_{2}-c_{4}\right) c_{6}-c_{3}c_{6}^{2}=0,  
\end{equation*}%
which yields  $c_{1}=c_{3}c_{5}$, $c_{2}=c_{4}$, $c_{6}=0.$ Hence we get from (5.18)%
\begin{equation}
gg^{\prime }=c_{7},\text{ }c_{7}\in \mathbb{R},\text{ }c_{7}\neq 0.  \tag{5.19}
\end{equation}
Substituting $\left( 5.19\right) $ into $\left( 5.4\right) $ yields%
\begin{equation}
g^{2}+\frac{1}{\left( f^{\prime }\right) ^{2}}+g^{2}\left( \frac{ff^{\prime
\prime }}{\left( f^{\prime }\right) ^{2}}-2\right) +\frac{2\lambda g^{2}}{%
c_{7}}\left( \frac{f^{\prime \prime }}{f\left( f^{\prime }\right) ^{2}}%
\right) +\frac{2\lambda g^{2}}{c_{7}f^{2}}=0.  \tag{5.20}
\end{equation}%
By dividing $\left( 5.20\right) $ with $g^{2}$ and after taking derivative
with respect to $z,$ we obtain the contradiction $g^{\prime }=0.$

\section{Some examples}

We illustrate some examples related with constant curvature factorable surfaces of type 2 in $\mathbb{I%
}^{3}.$

\begin{example}
Consider the factorable surfaces of type 2 in $\mathbb{I}^{3}$ given by

\begin{enumerate}
\item $\Phi _{3}:x=y\tan z,$ $\left( y,z\right) \in \left[ 0,\frac{\pi }{3}%
\right] ,$ (isotropic minimal),

\item $\Phi _{3}:x=-\sqrt{z},$ $\left( y,z\right) \in \left[ 0,2\pi \right] ,
$ $(H=-1)$,

\item $\Phi _{3}:x=-\frac{y^{2}}{4z},$ $\left( y,z\right) \in \left[ 1,1.4%
\right] \times \left[ 1,2\pi \right] ,$ (isotropic flat),

\item $\Phi _{3}:x=\frac{z}{y},$ $\left( y,z\right) \in \left[ 1,\pi \right]
\times \left[ 1,2\pi \right] ,$ $(K=-1).$
\end{enumerate}
The surfaces can be respectively plotted by Wolfram Mathematica 7.0 as in
Fig.1, ..., Fig.4.
\end{example}

\begin{figure}[ht]
\begin{center}
\includegraphics[scale=0.2]{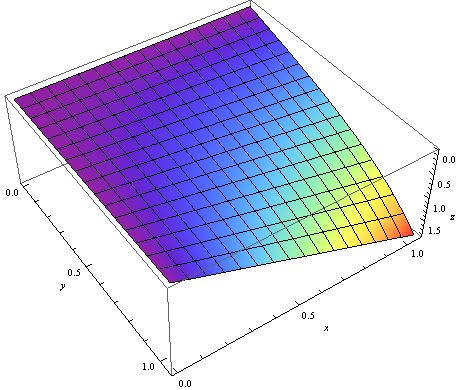}
\caption{An isotropic minimal factorable surface of type 2.}
\end{center}
\end{figure}

\begin{figure}[ht]
\begin{center}
\includegraphics[scale=0.2]{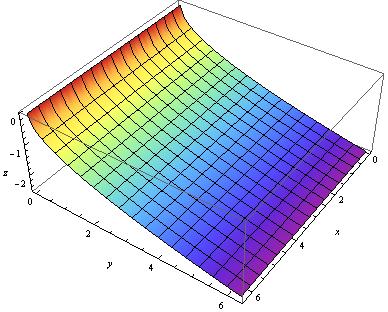}
\caption{A factorable surface of type 2 with $H=-1$.}
\end{center}
\end{figure}

\begin{figure}[ht]
\begin{center}
\includegraphics[scale=0.2]{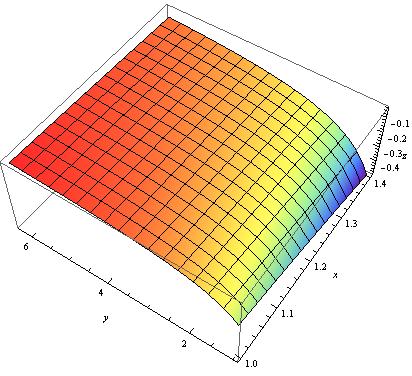}
\caption{An isotropic flat factorable surface of type 2.}
\end{center}
\end{figure}

\begin{figure}[ht]
\begin{center}
\includegraphics[scale=0.2]{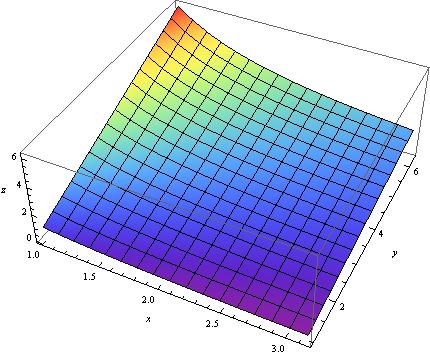}
\caption{A factorable surface of type 2 with $K=-1$.}
\end{center}
\end{figure}

\end{document}